\documentclass[makeidx]{amsart}
\usepackage{color}
\usepackage{amsmath}
\usepackage{amssymb}
\usepackage{graphicx}
\usepackage{hyperref}
\newtheorem{Proposition}{Proposition}[section]

  \newtheorem{Theorem}{Theorem}[section]
 
 \newtheorem{Definition}[Proposition]{Definition}
 \newtheorem{Note}[Proposition]{Note}

\newcommand {\z}{{\noindent}}
\def\blackslug{\hbox{\hskip 1pt \vrule width 4pt height 8pt depth 1.5pt
\hskip 1pt}}
\def\qed{\quad\blackslug\lower 8.5pt\null\par}

\def\CC{\mathbb{C}}
 \def\RR{\mathbb{R}}
 \def\NN{\mathbb{N}}

\def\Re{\mathrm{Re}}
\def\Im{\mathrm{Im}}

\makeindex
\begin{document}

\author{O. Costin$^1$ and X. Xia$^1$ }\title{From Taylor series of analytic functions to their global
  analysis}
\gdef\shortauthors{O.  Costin and X. Xia} \gdef\shorttitle{Global analytic
  reconstruction} \thanks{$1$.  Department of Mathematics, Ohio State
  University, 231 W. 18th Ave., Columbus, OH, 43210.}
  \maketitle

\bigskip

\begin{abstract}
 We analyze the conditions on the Taylor coefficients of an analytic function to admit global analytic continuation, complementing a recent paper of Breuer and Simon on general conditions for natural boundaries to form.  A new summation method is introduced to convert a relatively wide family of infinite sums and local expansions into integrals. The integral
  representations yield global information such as analytic
  continuability, position of singularities, asymptotics for large values of
  the variable and asymptotic location of zeros.

\end{abstract}

\setcounter{equation}{0}\section{Introduction}

Finding the global behavior of an analytic function in terms of its Taylor
coefficients is a notoriously difficult problem, in fact one which is impossible in full generality since undecidable statements can be formulated in these terms. However, a very interesting and quite general criterion for the disk of convergence of a Taylor series to coincide with its maximal domain of analyticity was recently discovered by Breuer and Simon \cite{Simon} (see also \cite{Sauzin1}).  The present paper complements this result by finding  criteria on the Taylor coefficients, say at zero, for the associated analytic function not to have natural boundaries and  to belong to the
class $\mathcal{M}$ of functions analytic in the complex plane with finitely
many cuts and with algebraic behavior at infinity (see Definition \ref{Def21}
below). Our condition  is that the coefficients $c_k$ admit generalized Borel summable (or Ecalle-Borel summable, EB) transseries in $k$. Many general classes of problems  in analysis are known to have EB transseries solutions.  For details on generalized Borel summability, transseries and resurgence see \cite{Book,Ecalle-book,Ecalle,Ecalle2,Ecalle3}.  

In particular  it is known \cite{Braaksma, Ecalle, Book} that, if  the $c_k$ are  solutions of generic linear or
nonlinear recurrence relations of arbitrary but finite order with analytic coefficients, then they are
are EB-summable.  Recurrence relations exist for instance
when the coefficients are obtained by solving differential equations by power
series. In a forthcoming paper we show that the Taylor coefficients of the Borel transform of solutions of generic systems of linear or nonlinear ODEs (in the setting of \cite{Duke}) also admit EB summable transseries. The global analytic structure of the Borel transform is crucial in understanding the monodromy of solutions of such equations.

 We also extend  our procedure to analyze the global behavior of entire functions, and to formal series, giving criteria directly  on the coefficients for the formal series to be Borel summable.

Globally reconstructing  function from its Taylor coefficients, when these admit EB summable transseries is effective, constructive and explicit -in the sense of producing integral representations far easier to analyze than the sums; this provides a new summation procedure, generalizing in some ways the Poisson summation formula.

 We recently used this approach to analyze a class of linear PDEs with variable coefficients,
\cite{CHT}.  One can obtain explicit integral representations for solutions of ODEs not known to be solvable such as
\begin{equation}
  \label{eq:thfilm}
  A\eta^2f^{(4)} +2A\eta f'''+\frac12 \eta f'-(1+a) f=0
\end{equation}
arising as the scaling pinching profile  $ h\sim (t_c-t)f(x(t_c-t)^{-1/2})$ of the thin film
equation, $$h_t+(h_{xxx} h)_x=0, \ h\sim (t_c-t)f(x(t_c-t)^{-1/2})$$ 
  where $t_c$ is the singularity time, where the interest is in the solution  {\em analytic at zero, $f_0$}. Let $A\neq 0$, and $a+\frac{1}{2} \notin \NN$. While it is not clear how to obtain representations of $f_0$ itself, the Taylor coefficients of $f_0$ are explicit. From this, our technique introduced in \S\ref{entire} by first considering the Laplace transform $F(p) = \mathcal{L}f (p)$, which is expressible in terms of integrals of Whittaker functions: with
  \begin{equation}\label{whit}
  F(p)= C e^{-\frac{p^{-2}}{8A}}p^{-3/2} \left[M_{-a-\frac{3}{4},\frac{1}{4}}\left(\frac{p^{-2}}{4A}\right)+ \frac{a\Gamma(a)}{2\sqrt{\pi}} W_{-a-\frac{3}{4},\frac{1}{4}}\left(\frac{p^{-2}}{4 A}\right) \right]
  \end{equation}
  we have $f = \mathcal{L}^{-1} F$. The proof of \eqref{whit} is sketched in \S\ref{4thoode}.

In particular, if $c_k=\varphi(k)$ where the function $\varphi$,  defined in the right
half plane,  is inverse Laplace transformable and its inverse Laplace
transform $\mathcal{L}^{-1} \varphi$ can be calculated in closed form, the function
$f$ has integral representations in terms of $\varphi$.  We will use some particularly simple examples for illustration. For the first one, the generalized Hurwitz zeta function,  our procedure quickly yields one of the known integral representations. For the other three, our procedure yields integral representations while  the global behavior of $\sum_k c_k z^k$   does not follow in any other obvious way: 
\begin{equation}
  \label{eq:ln2}
 c_k^{[1]}=\frac{1}{(k+a)^b};  \ \ c_k^{[2]} =\frac{1}{k^b+\ln
    k},\ \  c_k^{[3]}= \frac{1}{k^{k+1}}, \
  c_k^{[4]}=e^{\sqrt{k}}, \ (a,b>0)
\end{equation}
We find that
\begin{equation}
  \label{eq:sum0}
  f_1(z):=\sum_{k=1}^\infty c_k^{[1]}z^k= \frac{ z}{\Gamma(b)}\int_0^{\infty}\frac{[\ln(1+t)]^{b-1} \mathrm{d}t}{(1+t)^{a+1}(t-(z-1))}
\end{equation}
On the first Riemann sheet $f_1$ has only one singularity, at $z=1$, of
logarithmic type, and $f_1=o(z)$ for large $z$.  General Riemann surface
information and monodromy follow straightforwardly from \eqref{eq:sum0}.  A
similar complex analytic structure is shared by $f_2=\sum_{k=1}^\infty
c_k^{[2]}z^k$, which has one singularity at $z=1$ where it is analytic in
$\ln(1-z)$ and $(1-z)$; more precisely,
\begin{equation}
  \label{intlog}
 f_2(z)=-\frac{1}{2\pi i} \frac{\ln \ln
   z}{z}\oint_{0}^{\infty}\frac{e^{-u\ln(z)}}{(-u)^{b}+\ln (-u)}\mathrm{d}u +E(z)
\end{equation}
see Definition \ref{def2}, where $E$ is entire.

  The function $f_3(z)=\sum_{k=1}^\infty
c_k^{[3]}z^k$ is entire; questions answered regard
say the behavior for large negative $z$ (certainly not obvious from the series) or the asymptotic location of
zeros. It will follow that $f_3$ can be written as
\begin{equation}
  \label{eq:invl}
 f_3(z)=\int_0^{\infty}(1+u)^{-1}G(\ln(1+u))\left[\exp\left({\frac{ze^{-1}}{1+u}}\right)-1\right]\mathrm{d}u
\end{equation}
where $G(p)=s'_2(1+p)-s'_1(1+p)$ and $s_{1,2}$ are two branches of the
functional inverse of $s-\ln s$, cf. \S\ref{S115}.  Using the integral
representation of $f_{3}$, its behavior for large $z$ can be obtained from
(\ref{eq:invl}) by standard asymptotics methods; in particular, for large
negative $z$, $f_3$ behaves like a constant plus $z^{-1/2}e^{-z/e}$ times a
factorially divergent series (whose terms can be calculated).

For 
$ c_k^{[4]}$ we find 
\begin{equation}
  \label{ff4}
f_4(z)=\sum_{k=1}^\infty c_k^{[4]}z^k =-\frac{z}{2\sqrt{\pi}}\int_{C_1}\frac{p^{-3/2}e^{-\frac{1}{4p}}\mathrm{d}p}{e^p-z}
\end{equation}
where $C_1$ is a spiral $S_1$ followed by $[1,\infty)$, where $S_1$ starts at $0$ and ends at $1$, and is given  in polar coordinates by $r = \theta e^{2\pi i \theta}$, $\left(\theta \in [0,1]\right)$.

We also show that Borel summation of divergent series or transseries of
resurgent functions with finitely many Borel-plane singularities, as well as
the Abel-Plana version of the Euler-Maclaurin summation formula (see also
\cite{Stavros1}) can be derived by a natural extension of our analysis. Another illustration is obtaining the closed form Borel summed formula for $\ln \Gamma$, cf. \eqref{eq:lngamma} below.

A separate category is represented by lacunary series. Their coefficients do
not satisfy our assumption; however a slightly different approach allows for a
detailed study of the associated functions as the natural boundary is
approached, \cite{Advances}.

\section{Main results}

$ $\ \ \ {\em A first class of problems} is finding the location and type of
singularities in $\CC$ and the behavior for large values of the variable of
functions given by series with finite radius of convergence (Theorem
\ref{T1}), such as the first, second and fourth in (\ref{eq:ln2}).

{\em The second class of problems} amenable to the techniques presented
concerns the behavior at infinity (growth, decay, asymptotic location of
zeros etc.) of entire functions presented as Taylor series (Theorem \ref{T2}).

{\em The third class of problems} is essentially the converse of the two above: given a
function that has analytic continuation on some Riemann surface, how is this
reflected on $c_k$? (Theorem \ref{T1}.)

{\em The fourth class of problems} is to determine Borel summability of
series with zero radius of convergence such as
\begin{equation}
  \label{eq:entire}
  \tilde{f}_5=\sum_{n=0}^{\infty}n^{n+1} z^n
\end{equation}
in which
the coefficients of the series are analyzable (Theorem \ref{T3}).

\begin{Definition}\label{Def21}
  {\rm Let $\{a_j : 1\leq j \leq N\}$ be a set of nonzero complex numbers with distinct arguments. Let $\mathcal{M}$ consist of the functions algebraically bounded
at $\infty$ and  analytic in $\CC\setminus \bigcup_{j=1}^N
 \{a_j t: t\geq 1\}$ .  By dividing by a power of $z$ and subtracting out the principal part ({\em i.e.,}the negative powers of $z$) we can assume that
  $f\in\mathcal{M}'=\{f\in \mathcal{M}:f(z)=o(z)\ \text{as}\
  |z|\to\infty\}$. }\\
  
  {\rm This is one of the simplest settings often occurring in applications. We can see later from the proof that the approach is more general.}
\end{Definition}

\begin{Definition}
\label{def2}
{\rm Assume $g(s)$ is analytic in $U_\delta \backslash [0,\infty) $ for some $\delta > 0$, where $U_\delta = \{z : |\Im (z)| \leq \delta, \Re(z) \geq -\delta \}$ and $g(s) \to 0$ uniformly in $U_\delta$, as $\Re(s) \to \infty$.
Assume $\epsilon \leq \delta$. We define $\Gamma_\epsilon$ to be the contour around $\RR^+$ consisting of two rays $l_{1,\epsilon}, l_{2,\epsilon}$ and a semicircle $\gamma_{\epsilon}$, where $l_{1,\epsilon} = \{x - \epsilon i: a\in[0,\infty)\}$ oriented towards the left, $l_{2,\epsilon} = \{x + \epsilon i: a\in[0,\infty)\}$ oriented towards the right; $\gamma_{\epsilon}$ is the left semicircle centered at origin oriented clockwise.
Assume also that $g(s)$ is absolutely integrable over $\Gamma_\epsilon$ for some $\epsilon$.
We denote by
\begin{equation}
  \label{eq:sg}
  \oint_0^{\infty}g(s)\mathrm{d}s 
\end{equation}
the integral of $g$ over $\Gamma_\epsilon$. Since $g(s)$ vanishes at $\infty$, the integral is independent of the choice of $\epsilon$ as long as it is small enough.}
\end{Definition}

The following observations will simplify our proofs.

\begin{Note}\label{residue}
{\rm Let $g(s)$, $U_\delta$ and $\Gamma_\epsilon$ be as in Definition \ref{def2}. $\Gamma_\epsilon$ separates $\CC \backslash \Gamma_\epsilon$ into two regions. We denote the region containing $\RR^+$ by $S_1$ and the other by $S_2$. Let} 
\begin{align}
G_1(z) = \int_{\Gamma_\epsilon} \frac{g(s)\mathrm{d}s}{s-z} \quad\quad (z\in S_1)\\
G_2(z) = \int_{\Gamma_\epsilon} \frac{g(s)\mathrm{d}s}{s-z} \quad\quad (z\in S_2)
\end{align}
{\rm Then $G_1$ is analytic in $S_1$ and $G_2$ is analytic in $S_2$. By slightly deforming $\Gamma_\epsilon$ we are able to see that each $G_i$ can be analytically continued to $S_i \cup \Gamma_\epsilon$, $i = 1,2$. On $\Gamma_\epsilon$ their analytic continuations satisfy}
\begin{equation}
G_2(z)- G_1(z) = 2\pi i g(z)	\quad\quad (z \in \Gamma _\epsilon)
\end{equation}
{\rm Hence $G_2(z)$ can be analytically continued to $\CC \backslash [0,\infty)$ and $G_1(z)$ can be analytically continued to at least $U_\delta$ and in regions where $g$ is analytic. For each $z \in \CC \backslash [0,\infty)$}
\begin{equation}
G_2(z) = \oint_0^{\infty} \frac{g(s) \rm{d}s}{s-z} 
\end{equation}
\end{Note}

\begin{Note}
  {\rm A representation of the form \eqref{eq:sg} exists for Hilbert-transform-like integrals such as $h(t)= \int_0^{\infty}(s-t)^{-1}H(s)ds$ with $H$ analytic at
    zero, for instance  $h(t)=- (2\pi i)^{-1}
    \oint_0^{\infty}(s-t)^{-1}H(s)\ln s\,ds$.}
\end{Note}

\begin{Note}\label{COV}
  {\rm Consider the composition of $g$ with $s \mapsto \ln (1+s)$, the branch cut of which is chosen to be $(-\infty, -1]$. If $g$ is analytic in $U_\delta \setminus [0,\infty)$, then $g(\ln(1+s))$ is analytic in the set $-1+\exp(U_\delta \setminus [0,\infty))$. If in addition we have the decay condition $g(\ln(1+s)) = o(|s|^{-\alpha})$ for some $\alpha > 0$ as $|s| \to \infty$, then there exists a $\tilde{\delta}$ small enough such that $U_{\tilde{\delta}} \subseteq  -1+\exp(U_\delta \setminus [0,\infty))$. It is easy to see from the decay condition that}
  \begin{equation*}
  \int_{\Gamma_{\delta}} g(p) \mathrm{d}p = \int_{\exp(\Gamma_\delta) -1} \frac{g(\ln(1+s))}{1+s} \mathrm{d}s = \int_{\Gamma_{\tilde{\delta}}} \frac{g(\ln(1+s))}{1+s} \mathrm{d}s 
  \end{equation*}
 \rm{and thus we can make the change of variable}
 \begin{equation}
 \oint_{0}^{\infty} g(p) \mathrm{d}p = \oint_{0}^{\infty} \frac{g(\ln(1+s))}{1+s} \mathrm{d}s 
 \end{equation}
\end{Note}

While providing integral formulae in terms of functions with known
singularities which are often rather explicit, the following result can also
be interpreted as a {\em duality of resurgence}.  \footnote{After developing these methods,
  it has been brought to our attention that a duality between resurgent
  functions and resurgent Taylor coefficients has been noted in an unpublished
  manuscript by \'Ecalle. This will be further explored in a forthcoming paper.}. 
\begin{Theorem}\label{T1}
  (i) Assume that $f(z)=\sum_{k=0}^{\infty}c_k z^k$ is a series with positive, finite radius of convergence,  with $c_k$ having Borel sum-like representations of the form
\begin{equation}
  \label{eq:trs1}
  c_k=\sum_{j=1}^N a_j^{-k}\oint_0^{\infty}e^{-kp}F_j(p)\mathrm{d}p \hspace{0.4in} (k \ge 1)
\end{equation}
(\ref{eq:trs1}) with $a_j$ as in Definition 2.1, $F_j$ analytic in $U_\delta \backslash [0,\infty) $ for some $\delta > 0$ and algebraically bounded at $\infty$. Then, $f$ is given by
\begin{equation}
  \label{eq:recon3}
 f(z)=f(0)+z\oint_0^{\infty}\sum_{j=1}^N \frac{F_j(\ln (1+s))\mathrm{d}s}{(1+s)((1+s) a_j - z)}
\end{equation}

\z (ii) Furthermore, $f \in \mathcal{M}'$.The behavior of $f$ at $a_j$ and is of the same type as the behavior of $F_j(\ln(1+s))$ at $0$. More precisely, for small $z \notin [0,\infty)$,
\begin{equation}
\label{eq:singtype}
f(a_j (z+1)) = 2\pi i F_j (\ln (1+z)) + G(z)
\end{equation}
where $G(z)$ is analytic at 0.

\z (iii) Conversely, assume $f \in \mathcal{M}'$, and has finitely many singularities located at $\{a_j t_{j,l} \}$, $(1 \leq j \leq N, 1\leq l)$, with $1 =t _{j,1}$ and $ t _{j,l} < t_{j,l+1}$ for all $j, l$. Let
$c_k=f^{(k)}(0)/k!$; then $c_k$ have Borel sum-like representations of the form
\begin{equation}
  \label{eq:trs12}
  c_k=\frac{1}{2\pi i}\sum_{j=1}^N (a_j )^{-k}\oint_0^{\infty}e^{-ks}\ f( a_{j} e^s) \, \rm{d}s,\quad  k \geq 1
\end{equation}
\end{Theorem}

 The behavior at $a_j$ and
at $\infty$ will follow from the proof.

As it will be clear from the proofs, the method and results would apply, with
minor adaptations to functions of several complex variables.

\subsection{Entire functions}\label{entire}

 We
restrict the analysis to entire functions of exponential order one, with
complete information on the Taylor coefficients. Such functions include of
course the exponential itself, or expressions such as $f_3$. It is useful to
start with $f_3$ as an example. The analysis is brought to the case in
Theorem \ref{T1} by first taking a Laplace transform. Note that
\begin{equation}
  \label{eq:eq41}
  \int_0^{\infty}e^{-xz}f(z)\mathrm{d}z=\frac{1}{x}\sum_{n=1}^\infty\frac{n!}{n^{n+1}x^{n}}
\end{equation}
The study of entire functions of exponential order one likely involves the
factorial, and then a Borel summed representation of the Stirling formula is
needed; this is provided in the Appendix.
\begin{Theorem}\label{T2}
  Assume that the entire function $f$ is given by
\begin{equation}
  \label{eq:enti22}
 f(z)= \sum_{k=1}^\infty \frac{c_k z^k }{k!}
\end{equation}
with $c_k$ as in Theorem~\ref{T1} (i). Then,
\begin{equation}
  \label{eq:recon31}
 f(z)=
  \oint_0^{\infty}\sum_{j=1}^N \left[\left(e^{\frac{z}{a_j(1+s)}}-1\right)
\frac{F_j(\ln (1+s))}{(1+s)}\right]\mathrm{d}s
\end{equation}
\end{Theorem}
As in the simple example, the behavior at infinity follows from the integral
representation by classical means.
\subsection{Borel summation}
We obtain from Theorem~\ref{T1}, in the same way as above, the
following.
\begin{Theorem}\label{T3}
  Consider the formal power series
\begin{equation}
  \label{eq:enti23}
 \tilde{f}(z)= \sum_{k=1}^\infty {c_k k! x^{-k-1}}
\end{equation}
with coefficients $c_k$ as in Theorem~\ref{T1} (i).
Then the series \eqref{eq:enti23} is (generalized)  Borel summable to
\begin{multline}
  \label{eq:bs1}
\int_0^{\infty}dp e^{-px}p\sum_{j=1}^N\oint_0^{\infty}\frac{F_j(\ln(1+s))}{(1+s)(a_js+a_j-p)}\mathrm{d}s\\=
   \sum_{j=1}^N\oint_0^{\infty}\frac{F_j(\ln (1+s))}{1+s}\left(-\frac{1}{x}+a_j(s+1)e^{-a_j(s+1)x}\,\mathrm{Ei}\left(a_j(s+1){x}\right)\right)\mathrm{d}s
\end{multline}
\end{Theorem}
The proof proceeds as in the previous sections, taking now a Borel
transform followed by Laplace transform.
\subsection{Other applications; the examples in the introduction}\label{Examples}
\subsubsection{Other growth rates}\label{Growthrates}

Series with coefficients with growth rates precluding a straightforward inverse Laplace transform can be accommodated, for instance by analytic continuation. 
We have for positive $\gamma$,
\begin{equation}
  \label{eq:eqinvl}
  e^{-\gamma\sqrt{
      n}}=\frac{\gamma}{2\sqrt{\pi}}\int_0^{\infty}p^{-3/2}e^{-\frac{\gamma^2}{4p}}e^{-np}\mathrm{d}p
\end{equation}
which can be analytically continued in $\gamma$. We note first that the
contour cannot be, for this function, detached from zero. Instead, we keep $[1,\infty)$  as part of the original contour fixed and, deform the part $[0,1]$ by simultaneously rotating $\gamma$ and $p$  to maintain $\gamma^2/p$ real
and positive near the
origin. We get 
\begin{equation}
  \label{eq:eqerfc}
  e^{\sqrt{n}}=-\frac{1}{2\sqrt{\pi}}\int_{C_1} p^{-3/2} e^{-\frac{1}{4p}}e^{-np}\mathrm{d}p
\end{equation}
 and \eqref{ff4} follows, for the same reason Theorem \ref{T1} (i) holds.  In particular, 
\begin{equation}
  \label{eq:limit1}
  \lim_{z\to -1^+}\sum_{n=1}^{\infty}e^{\sqrt{n}}z^n=
  \frac{1}{2\sqrt{\pi}}\int_{C_1}\frac{p^{-3/2}e^{-\frac{1}{4p}}}{e^p+1}\mathrm{d}p
\end{equation}
 The sum \eqref{eq:limit1} is  unwieldy numerically,
while the integral \eqref{eq:eqinvl} can be evaluated accurately by standard
means. In a similar way we get
\begin{equation}
  \label{eq:eqsum}
  \sum_{k=0}^{\infty}\frac{e^{i\sqrt{k}}}{k^{a}}=-\frac{\gamma
    2^{a-1/2}}{\sqrt{\pi}}
\int_C \frac{e^{-\frac{1}{8p}}U(2a+1/2;\frac{1}{\sqrt{2p}})}{p^{a-1}(e^p+1)}\mathrm{d}p
\end{equation}
for $a>1/2$ for which the series converges.
Here $C$ is a contour consisting of $S_2$ followed by $[1,\infty)$, where $S_2$ starts at $0$ and ends at $1$, and is given by $r = 1 - \theta/\pi$, $\left(\theta \in [0,\pi]\right)$ in polar coordinates. $U$ is the parabolic cylinder function
\cite{Abramowitz}.

{\em The coefficients $c_k^{[1]}$ in \eqref{eq:ln2}}.
 We have
\begin{equation}
  \label{eq:invlp3}
  \mathcal{L}^{-1}\left[\frac{1}{(n+a)^b}\right]=\Gamma(b)^{-1}p^{b-1}e^{-ap}
\end{equation}
The rest follows in the same way \eqref{ff4} was obtained, after changing
variables to $1+t=e^p$.

{\em The coefficients $c_k^{[2]}$}.
 We let $x=n$ and take the inverse
Laplace transform in $x$:
\begin{equation}
  \label{eq:dual1}
  G(p) = \frac{1}{2\pi i}\int_{1-i\infty}^{1+i\infty}\frac{e^{xp}}{x^{b}+\ln x}\mathrm{d}x
\end{equation}
where the contour can be bent backwards for $p \in \RR^+$, to hang around $\RR^-$. Then, with
the change of variable $x=-u$  (\ref{eq:dual1}) becomes
\begin{equation}
  \label{eq:dual2}
 G(p)= \frac{1}{2 \pi i}\oint_{0}^{\infty}\frac{e^{-up}}{(-u)^{b}+\ln (-u)}\mathrm{d}u
\end{equation}
and thus
\begin{equation*}
c_k = (\mathcal{L} G) (k) = \int_0^{\infty} G (p) e^{-kp}{\mathrm{d}p }  = \oint_0^{\infty} \left[ \frac{-G(p) \ln p}{ 2 \pi i} \right] e^{-kp}{\mathrm{d}p}
\end{equation*}
We see that $F_1(p) = (- G(p) \ln p)/{2 \pi i}$ and by Theorem 2.1
\begin{equation*}
f_2 (z) = z \oint_0^{\infty} \frac {\tilde{G}(s)}{s - (z-1)} {\mathrm{d}s}
\end{equation*}
where $\tilde {G}(s) = F_1 (\ln(1+s))/(1+s)$.
Hence the singularity of $f_2(z)$, at $z=1$ on the first Riemann sheet, according to Note \ref{residue} is
that of $ \phi(z) =  2 \pi i \tilde{G}(z - 1)$, as in
(\ref{intlog}).

The example of {\em the coefficients $c_k^{[3]}$} is studied in a similar way as Theorem
\ref{T2}; related calculations can be found in \S\ref{S115}.

{\em The coefficients $c_k^{[4]}$} were treated at the beginning of this
  section.

{\em Another example} is provided by the log of the Gamma function,
$\ln\Gamma(n)=\sum_{k=1}^n \ln k$. It is convenient to first subtract out the leading
behavior of the sum to arrange that the summand is inverse Laplace
transformable. With 
$$g_n=\ln\Gamma(n+1)-\left((n+1)\ln (n+1)-n-\frac{1}{2}\ln (n+1)\right)$$  we get 
\begin{multline}
  g_N=\sum_{1}^N
  \Big[1-\Big(\frac{1}{2}+n\Big)\ln\Big(1+\frac{1}{n}\Big)\Big]=\sum_{1}^N\int_0^\infty
  e^{-np} \frac{ 1-\frac{p}{2}-(\frac{p}{2}+1)e^{-p}}{p^2}\mathrm{d}p
\end{multline}
where $\mathcal{L}^{-1}$ of the summand in the middle term is most easily
obtained by noting that its second derivative is a rational function. 
Summing as usual $e^{-np}$ we get
\begin{equation}
  \label{eq:lngamma}
  \ln\Gamma(n)=n(\ln n -1)-\frac{1}{2}\ln n +\frac{1}{2}\ln(2\pi)+\int_0^{\infty}\frac{\displaystyle
  1-\frac{p}{2}-\Big(\frac{p}{2}+1\Big)e^{-p}}{p^2(e^{-p}-1)}e^{-np}\mathrm{d}p
\end{equation}

Obviously, if the behavior of the coefficients is of the form $A^kc_k$ where
$c_k$ satisfies
the conditions in the paper, one simply changes the independent variable to  $z'=Az$.

\section{Proof of Theorem \ref{T1}}
\z If $f\in\mathcal{M}'$ we write the Taylor coefficients in the form
\begin{equation}
      \label{eq:cf1}
    c_k=\frac{1}{2\pi i}\oint \frac{ f(p)\mathrm{d}p}{p^{k+1}}  \quad\quad(k\geq 1)
  \end{equation}
  where the contour of integration is a small circle of radius $r$ around the
  origin. We attempt to increase $r$ without bound. In the process, the contour
  will hang around the singularities of $f$ as shown in Figure 1. Each integral over a curve that wraps around a ray $\{a_j t:t\geq 1\}$ converges by the decay assumptions and the contribution of the arcs
  at large $r$ vanishes, since $f(z) = o(z)$ as $z \to \infty$.
  
\begin{figure}
  \includegraphics[width=12cm]{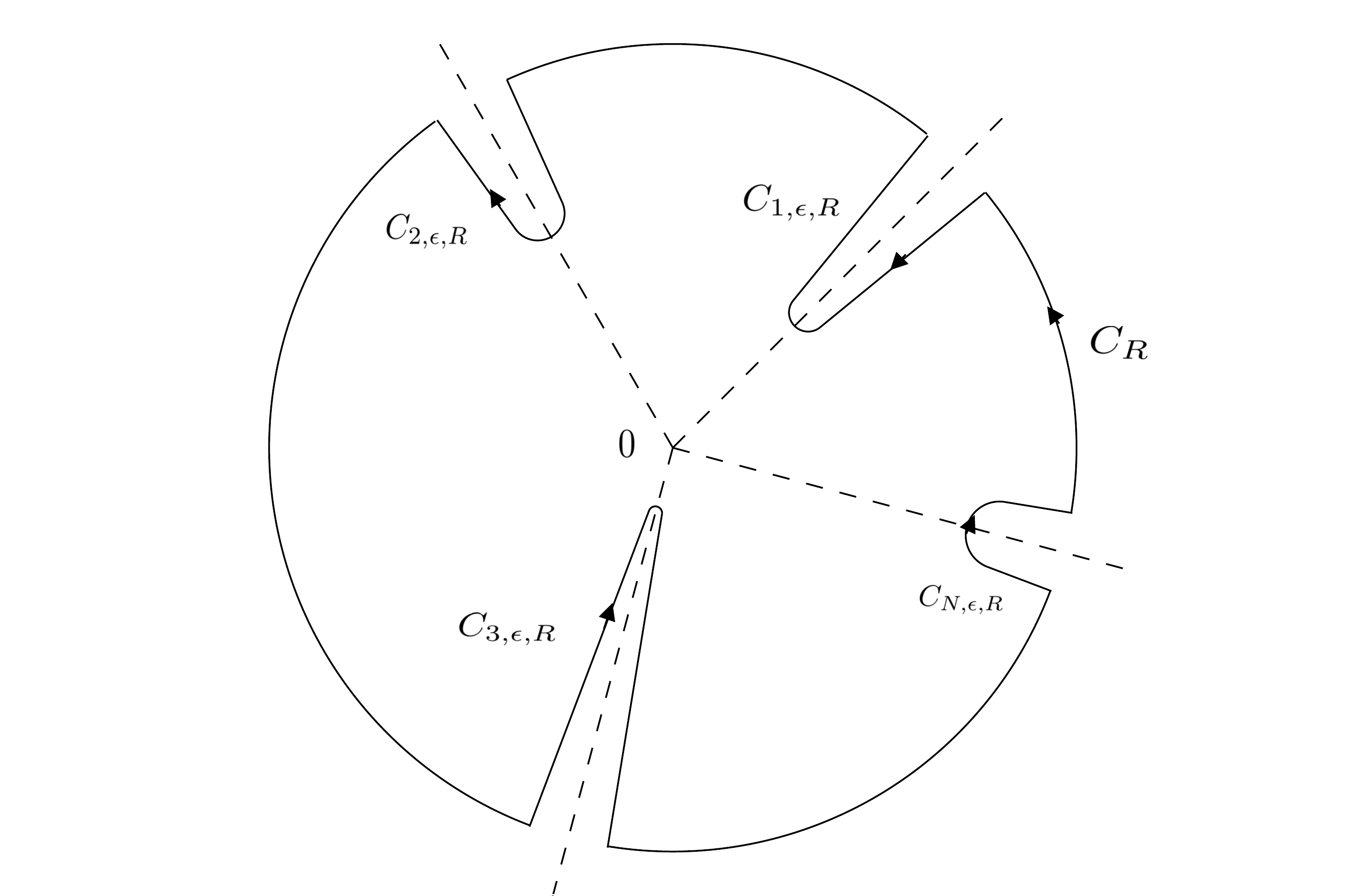}
  \caption{Singularities of $f$, cuts, direction of integration and Cauchy
  contour deformation. }
\end{figure}

  To be more precise, let $C_{j,\epsilon}$ be the part of  the image of $\Gamma_{\epsilon}$ under the mapping $s \to a_j e^s $, let $C_{j,\epsilon, R}$ be the part of $C_{j,\epsilon}$ inside the disk $|s| \leq R$, and $C_R$ be the part of the contour on $|s| = R$. Then for $\epsilon$ small enough and $R$ large enough we have
\begin{equation*}
c_k=\frac{1}{2\pi i}\oint \frac{ f(p)\mathrm{d}p}{p^{k+1}} = \frac{1}{2\pi i} \left( \sum_{j = 1}^{N} \int_{C_{j, \epsilon, R}} \frac{ f(p)\mathrm{d}p}{p^{k+1} }+ \int_{C_R} \frac{ f(p)\mathrm{d}p}{p^{k+1} }\right)
\end{equation*}
By the change of variable $p = a_j e^s$ and letting $R \to \infty$ we get
\begin{equation*}
\int_{a_j C_{\epsilon, R}} \frac{ f(p)\mathrm{d}p}{p^{k+1} } = \int_{\Gamma_{\epsilon}} \frac{a_j e^s  f(a_j e^s) \mathrm{d}s}{(a_j e^s)^{k+1} } = \oint_{0}^{\infty} {a_j}^{-k} e^{-ks} f(a_j e^s)  \mathrm{d}s
\end{equation*}  
and the integral over $C_R$ vanishes as $R \to \infty$ by decay condition for $k \geq 1$.

\bigskip

  \z In the opposite direction, first let $\epsilon$ be small enough so that for all $j$ and $k = 1$
  \begin{equation}
  \label{eq:path}
  \oint_{0}^{\infty}{e^{-kp}F_j(p) \mathrm{d}p} = \int_{\Gamma_{\epsilon}}{e^{-kp}F_j(p)  \mathrm{d}p} 
  \end{equation}
 Then for all $1\leq j \leq N$, $k \geq 1$ (\ref{eq:path}) is true. Also let $z$ be small so that 
  \begin{equation}
  \label{eq:est1}
  |{a_j}^{-1}e^{-p} z| \leq \delta^2 < 1
  \end{equation}
  for all $j$ and $p \in \Gamma_\epsilon$.

Then, by the dominated convergence theorem (which applies in this case, see (\ref{eq:estimate})) we have
\begin{align}
  \label{eq:converse}
  f(z) - f(0)&
  = \sum_{k=1}^{\infty}c_k z^k  
  =\sum_{k=1}^{\infty}\left({\sum_{j=1}^N a_j^{-k}\int_{\Gamma_\epsilon}e^{-kp}F_j(p)\mathrm{d}p} \right) z^k\\ \notag
&= \int_{\Gamma_\epsilon} \sum_{j=1}^{N} \left( \sum_{k=1}^{\infty}  (a_j^{-1}e^{-p} z)^k \right) F_j(p) \mathrm{d}p = \int_{\Gamma_\epsilon} \sum_{j=1}^{N} \left( \frac{a_j^{-1}e^{-p} z}{1 - a_j^{-1}e^{-p} z} \right) F_j(p) \mathrm{d}p\\ \notag
&=\sum_{j=1}^{N} \oint_0^{\infty} \frac{a_j^{-1}e^{-p} z}{1 - a_j^{-1}e^{-p} z} F_j(p) \mathrm{d}p =z\sum_{j=1}^N 
  \oint_0^{\infty}\frac{F_j(\ln (1+s))\mathrm{d}s}{(1+s)(s a_j+a_j-z)}
\end{align}
as stated.  The third equality holds because we have, in view of \eqref{eq:est1},
\begin{align}
\label{eq:estimate}
\left| \sum_{j=1}^N (a_j^{-1}e^{-p} z)^kF_j(p)\right| \leq \sum_{j=1}^N |a_j^{-1} e^{-p} z|^{k/2} \left( |a_j^{-1} e^{-p} z|^{k/2} |F_j(p)|\right) \nonumber\\
\leq \sum_{j=1}^N \delta^k \left( |a_j^{-1} e^{-p} z|^{k/2} |F_j(p)|\right)
\end{align}
For each $j$, $|a_j^{-1} e^{-p} z|^{k/2} F_j(p)$ is integrable over $\Gamma_\epsilon$ since $F_j$ is algebraically bounded at $\infty$, so we may interchange the order of integration and summation over $k$. The last equality holds because $F_j(\ln (1+s))/(s a_j+a_j-z) = o(|s|^{-\alpha})$ for some $\alpha > 1$ so we can make the change of variable $p=\ln (1+s)$, see Note \ref{COV}. Hence \eqref{eq:trs1} holds for $z$ small.\\
\indent Given $j\in \{1, ..., N\}$ we may write
\begin{equation}\label{intj}
I_j(z) := \oint_0^{\infty}\frac{F_j\left(\ln (1+s)\right) \rm{d}s}{(1+s)(s a_j+a_j-z)} = \oint_0^{\infty}\frac{F_j\left(\ln (1+s)\right)/\left(a_j (1+s)\right)\, \rm{d}s}{\left(s-\left(z/a_j - 1\right)\right)} 
\end{equation}
Since $F_j$ is analytic in $\CC \backslash [0,\infty)$ and is algebraically bounded at $\infty$, $F_j(\ln (1+s))/\left(a_j (1+s)\right)$ is analytic in $\CC \backslash [0,\infty)$ and vanishes uniformly as $\Re(s) \to \infty$. Then it becomes obvious from Note \ref{residue} that the integral $I_j$ in \eqref{intj} is analytic in $\CC \backslash \{a_j t: t \geq 1\}$. Thus $f(z)$ can be analytically continued to $\CC \backslash \bigcup_{j=1}^{N} \{a_j t: t \geq 1\}$.\\

To see that $f(z) = o(z)$ as $|z| \to \infty$, it suffices to show that for each $j$, $I_j(z) = o(1)$. Assume $3\epsilon < \delta$. We use the contours $\exp(\Gamma_\epsilon)-1$ and $\exp(\Gamma_{3{\epsilon}})-1$. If $|\arg(z/a_j)| \geq 2 {\epsilon}$ and $|z/a_j|$ is large enough we write:
\begin{equation}
I_j(z) = \int_{\exp(\Gamma_\epsilon) -1}\frac{F_j\left(\ln (1+s)\right)/\left(a_j (1+s)\right)\, \rm{d}s}{\left(s-\left(z/a_j - 1\right)\right)}
\end{equation}
Then there exist some positive number $\rho_1$ such that  $|s-\left(z/a_j - 1\right)| \geq \rho_1 |s|$ for each $s\in \exp(\Gamma_\epsilon) -1$, so we can use dominated convergence  to obtain $I_j(z) \to 0$ as $|z| \to \infty$. If $|\arg(z/a_j)| \leq 2 {\epsilon}$ we use Note \ref{residue} to write:
\begin{equation}\label{insideIj}
I_j(z) = \int_{\Gamma_{3{\epsilon}}}\frac{F_j\left(\ln (1+s)\right)/\left(a_j (1+s)\right)\, \rm{d}s}{\left(s-\left(z/a_j - 1\right)\right)} + 2\pi i F_j(\ln(z/a_j))/z
\end{equation}
Then there exist some positive number $\rho_2$ such that  $|s-\left(z/a_j - 1\right)| \geq \rho_2 |s|$ for each $s\in \exp(\Gamma_3\epsilon) -1$, so we can use dominated convergence to prove the integral in the right hand side of \eqref{insideIj} is $o(1)$. By assumption it is obvious that $F_j(\ln(z/a_j))/z = o(1)$, so in this case we also have $I_j(z) \to 0$ as $|z| \to \infty$.\\

The nature of the singularities of $f$ is derived from Note \ref{residue}. Let $z\notin [0,\infty)$ be small, then for each $j \in \{1,...,N\}$
\begin{align*}
&f((z+1)a_j)\\ &=f(0) + ((z+1)a_j)\sum_{l \neq j} I_l((z+1)a_j)+ ((z+1)a_j)I_j((z+1)a_j)\\
				   &=f(0) + ((z+1)a_j)\sum_{l \neq j} I_l((z+1)a_j)\\
				    &\quad\quad+((z+1)a_j)\oint_0^{\infty}\frac{F_j\left(\ln (1+s)\right)/\left(a_j (1+s)\right)\, \rm{d}s}{s-z}\\
				    &=f(0) + ((z+1)a_j)\sum_{l \neq j} I_l((z+1)a_j) + A(t) + 2 \pi i F_j(\ln(1+z))\end{align*}
where $A(t)$ is analytic at $z=0$. The last equality is obtained by Note \ref{residue}. It is obvious from Note \ref{residue} that each $I_l((z+1)a_j)$ $(l\neq j)$ is analytic on $[0, \infty)$. Thus 
\begin{equation}
f((z+1)a_j) = \tilde{G}(z) + 2 \pi i F_j(\ln(1+z))
\end{equation}
where $\tilde{G}(z)$ is analytic at $z=0$. Hence \eqref{eq:singtype} follows.

\section{Appendix}
\label{EB}

\subsection{Simple integral representation of  $1/n!$ \cite{Abramowitz}} 
We have
\begin{equation}
  \label{eq:eqg}
  \frac{1}{\Gamma(z)}=-\frac{ie^{-\pi i
      z}}{2\pi}\oint_0^{\infty}s^{-z}e^{-s}ds=-\frac{ie^{-\pi i
      z}z^{-z}}{2\pi}\oint_0^{\infty}s^{-z}e^{-zs}ds
\end{equation}
with our convention of contour integration. From here, one can proceed 
as in \S\ref{S115}.

\subsection{The Gamma function and Borel summed Stirling formula}\label{S115}
We have
\begin{multline}
  n!=\int_0^\infty t^ne^{-t}\mathrm{d}t=n^{n+1}\int_0^\infty e^{-n(s-\ln s)}\mathrm{d}s\\=
  n^{n+1}\int_0^1 e^{-n(s-\ln s)}\mathrm{d}s+ n^{n+1}\int_1^\infty e^{-n(s-\ln s)}\mathrm{d}s
\end{multline}

\z On $(0,1)$ and $(1,\infty)$ separately, the function $s-\ln(s)$ is
monotonic and we may write, after inverting $s-\ln(s)=t$ on the two
intervals to get $s_{1,2}=s_{1,2}(t)$ \footnote{The functions $s_{1,2}$ are
  given by branches of $-W(-e^t)$, where $W$ is the Lambert function.},
\begin{equation}\label{e25}
 n!= n^{n+1}\int_1^{\infty}e^{-nt}(s'_2(t)-s'_1(t))\mathrm{d}t=n^{n+1}e^{-n}
\int_0^{\infty}e^{-np}G(p)\mathrm{d}p
\end{equation}
\z where $G(p)=s_2'(1+p)-s_1'(1+p)$. From the definition it follows that $G$
is bounded at infinity and $p^{1/2}G$ is analytic in $p$ at $p=0$.
Using now (\ref{e25}) and Theorem~\ref{T1} in (\ref{eq:eq41}) we get
\begin{equation}
  \label{eq:lapl41}
  \int_0^{\infty}e^{-xz}f_3(z)\mathrm{d}z=\frac{1}{ x^2}\int_0^{\infty}\frac{G(\ln(1+t))}{(te+(e-x^{-1}))(t+1)}\mathrm{d}t
\end{equation}
Upon taking the inverse Laplace transform we obtain (\ref{eq:invl}).

\subsection{Solution of \eqref{eq:thfilm}}\label{4thoode}
Assume the solution to \eqref{eq:thfilm} which is analytic at $\eta=0$ has Taylor expansion
$f(\eta) = \sum_{k=0}^{\infty} c_k \eta^k$. Then 
\begin{equation}
c_k = \frac{c_1 (-1)^{k/2-1/2}}{A\Gamma(-1/2-a)} \frac{\Gamma(k/2-1-a) k}{\left(\Gamma(k+1)\right)^2} \quad\quad (k\,\,\rm{is}\,\,odd)
\end{equation}
\begin{equation}
c_k = 0 \quad\quad (k\,\,\rm{is}\,\,even)
\end{equation}
It is obvious that $f(\eta)$ is entire. Consider the Laplace transform $F(p) = \mathcal{L}f(p)$. Then $F(p) = \sum_{k=0}^{\infty} \Gamma(k+1) c_k p^{-k-1}$. Let $G(p) = F(1/p)$; then $G(z)$ is a solution to the differential equation
\begin{equation}
G''(z)+\left(-\frac{4}{z}+\frac{z}{2A}\right) G'(z) + \left(\frac{6}{z^2}-\frac{3/2+a}{A}\right)G(z) = 0
\end{equation}
and $G(z)$ is analytic at $z=0$. 

The normalization transformation  $G(z)=z^{3/2}e^{-z^2/A}h(z^2/A)$  (cf. \cite{Duke} for a general description of normalization methods) yields
  \begin{equation}
    \label{eq:eqW}
    h''-\frac{7}{4}h'+\frac{1}{16}\left(12-\frac{3+4a}{s}+\frac{3}{s^2}\right)h=0
  \end{equation}
solvable in terms of Whittaker functions \cite{Abramowitz}.  Substituting back, by straightforward algebra, this yields \eqref{whit}.

\subsection{Acknowledgments}  M. Kruskal introduced OC to questions of the
type addressed in the paper; OC is grateful to  R. D. Costin, J. McNeal, D. Sauzin,  and  S. Tanveer for very useful discussions. Work supported in part
by NSF grants    DMS-0600369 and DMS-1108794.

\end{document}